\documentclass[reqno,12pt,a4paper]{amsart}

\usepackage{mathrsfs}
\usepackage{amssymb}
\usepackage{amsfonts}
\usepackage{latexsym}
\usepackage{amsthm}

\usepackage{graphicx}
\def\lb{\label}

\newcommand{\er}[1]{\textrm{(\ref{#1})}}

\begin{document}


\renewcommand{\theequation}{\arabic{section}.\arabic{equation}}
\theoremstyle{plain}
\newtheorem{theorem}{\bf Theorem}[section]
\newtheorem{lemma}[theorem]{\bf Lemma}
\newtheorem{corollary}[theorem]{\bf Corollary}
\newtheorem{proposition}[theorem]{\bf Proposition}
\newtheorem{definition}[theorem]{\bf Definition}
\newtheorem{remark}[theorem]{\it Remark}

\def\a{\alpha}  \def\cA{{\mathcal A}}     \def\bA{{\bf A}}  \def\mA{{\mathscr A}}
\def\b{\beta}   \def\cB{{\mathcal B}}     \def\bB{{\bf B}}  \def\mB{{\mathscr B}}
\def\g{\gamma}  \def\cC{{\mathcal C}}     \def\bC{{\bf C}}  \def\mC{{\mathscr C}}
\def\G{\Gamma}  \def\cD{{\mathcal D}}     \def\bD{{\bf D}}  \def\mD{{\mathscr D}}
\def\d{\delta}  \def\cE{{\mathcal E}}     \def\bE{{\bf E}}  \def\mE{{\mathscr E}}
\def\D{\Delta}  \def\cF{{\mathcal F}}     \def\bF{{\bf F}}  \def\mF{{\mathscr F}}
\def\c{\chi}    \def\cG{{\mathcal G}}     \def\bG{{\bf G}}  \def\mG{{\mathscr G}}
\def\z{\zeta}   \def\cH{{\mathcal H}}     \def\bH{{\bf H}}  \def\mH{{\mathscr H}}
\def\e{\eta}    \def\cI{{\mathcal I}}     \def\bI{{\bf I}}  \def\mI{{\mathscr I}}
\def\p{\psi}    \def\cJ{{\mathcal J}}     \def\bJ{{\bf J}}  \def\mJ{{\mathscr J}}
\def\vT{\Theta} \def\cK{{\mathcal K}}     \def\bK{{\bf K}}  \def\mK{{\mathscr K}}
\def\k{\kappa}  \def\cL{{\mathcal L}}     \def\bL{{\bf L}}  \def\mL{{\mathscr L}}
\def\l{\lambda} \def\cM{{\mathcal M}}     \def\bM{{\bf M}}  \def\mM{{\mathscr M}}
\def\L{\Lambda} \def\cN{{\mathcal N}}     \def\bN{{\bf N}}  \def\mN{{\mathscr N}}
\def\m{\mu}     \def\cO{{\mathcal O}}     \def\bO{{\bf O}}  \def\mO{{\mathscr O}}
\def\n{\nu}     \def\cP{{\mathcal P}}     \def\bP{{\bf P}}  \def\mP{{\mathscr P}}
\def\r{\rho}    \def\cQ{{\mathcal Q}}     \def\bQ{{\bf Q}}  \def\mQ{{\mathscr Q}}
\def\s{\sigma}  \def\cR{{\mathcal R}}     \def\bR{{\bf R}}  \def\mR{{\mathscr R}}
                \def\cS{{\mathcal S}}     \def\bS{{\bf S}}  \def\mS{{\mathscr S}}
\def\t{\tau}    \def\cT{{\mathcal T}}     \def\bT{{\bf T}}  \def\mT{{\mathscr T}}
\def\f{\phi}    \def\cU{{\mathcal U}}     \def\bU{{\bf U}}  \def\mU{{\mathscr U}}
\def\F{\Phi}    \def\cV{{\mathcal V}}     \def\bV{{\bf V}}  \def\mV{{\mathscr V}}
\def\P{\Psi}    \def\cW{{\mathcal W}}     \def\bW{{\bf W}}  \def\mW{{\mathscr W}}
\def\o{\omega}  \def\cX{{\mathcal X}}     \def\bX{{\bf X}}  \def\mX{{\mathscr X}}
\def\x{\xi}     \def\cY{{\mathcal Y}}     \def\bY{{\bf Y}}  \def\mY{{\mathscr Y}}
\def\X{\Xi}     \def\cZ{{\mathcal Z}}     \def\bZ{{\bf Z}}  \def\mZ{{\mathscr Z}}
\def\O{\Omega}

\newcommand{\gA}{\mathfrak{A}}
\newcommand{\gB}{\mathfrak{B}}
\newcommand{\gC}{\mathfrak{C}}
\newcommand{\gD}{\mathfrak{D}}
\newcommand{\gE}{\mathfrak{E}}
\newcommand{\gF}{\mathfrak{F}}
\newcommand{\gG}{\mathfrak{G}}
\newcommand{\gH}{\mathfrak{H}}
\newcommand{\gI}{\mathfrak{I}}
\newcommand{\gJ}{\mathfrak{J}}
\newcommand{\gK}{\mathfrak{K}}
\newcommand{\gL}{\mathfrak{L}}
\newcommand{\gM}{\mathfrak{M}}
\newcommand{\gN}{\mathfrak{N}}
\newcommand{\gO}{\mathfrak{O}}
\newcommand{\gP}{\mathfrak{P}}
\newcommand{\gQ}{\mathfrak{Q}}
\newcommand{\gR}{\mathfrak{R}}
\newcommand{\gS}{\mathfrak{S}}
\newcommand{\gT}{\mathfrak{T}}
\newcommand{\gU}{\mathfrak{U}}
\newcommand{\gV}{\mathfrak{V}}
\newcommand{\gW}{\mathfrak{W}}
\newcommand{\gX}{\mathfrak{X}}
\newcommand{\gY}{\mathfrak{Y}}
\newcommand{\gZ}{\mathfrak{Z}}

\def\ve{\varepsilon} \def\vt{\vartheta} \def\vp{\varphi}  \def\vk{\varkappa}

\def\Z{{\mathbb Z}} \def\R{{\mathbb R}} \def\C{{\mathbb C}}  \def\K{{\mathbb K}}
\def\T{{\mathbb T}} \def\N{{\mathbb N}} \def\dD{{\mathbb D}} \def\S{{\mathbb S}}
\def\B{{\mathbb B}}


\def\la{\leftarrow}              \def\ra{\rightarrow}     \def\Ra{\Rightarrow}
\def\ua{\uparrow}                \def\da{\downarrow}
\def\lra{\leftrightarrow}        \def\Lra{\Leftrightarrow}
\newcommand{\abs}[1]{\lvert#1\rvert}
\newcommand{\br}[1]{\left(#1\right)}

\def\lan{\langle} \def\ran{\rangle}


\def\lt{\biggl}                  \def\rt{\biggr}
\def\ol{\overline}               \def\wt{\widetilde}
\def\no{\noindent}


\let\ge\geqslant                 \let\le\leqslant
\def\lan{\langle}                \def\ran{\rangle}
\def\/{\over}                    \def\iy{\infty}
\def\sm{\setminus}               \def\es{\emptyset}
\def\ss{\subset}                 \def\ts{\times}
\def\pa{\partial}                \def\os{\oplus}
\def\om{\ominus}                 \def\ev{\equiv}
\def\iint{\int\!\!\!\int}        \def\iintt{\mathop{\int\!\!\int\!\!\dots\!\!\int}\limits}
\def\el2{\ell^{\,2}}             \def\1{1\!\!1}
\def\sh{\sharp}
\def\wh{\widehat}
\def\bs{\backslash}
\def\na{\nabla}

\def\sh{\mathop{\mathrm{sh}}\nolimits}
\def\all{\mathop{\mathrm{all}}\nolimits}
\def\Area{\mathop{\mathrm{Area}}\nolimits}
\def\arg{\mathop{\mathrm{arg}}\nolimits}
\def\const{\mathop{\mathrm{const}}\nolimits}
\def\det{\mathop{\mathrm{det}}\nolimits}
\def\diag{\mathop{\mathrm{diag}}\nolimits}
\def\diam{\mathop{\mathrm{diam}}\nolimits}
\def\dim{\mathop{\mathrm{dim}}\nolimits}
\def\dist{\mathop{\mathrm{dist}}\nolimits}
\def\Im{\mathop{\mathrm{Im}}\nolimits}
\def\Iso{\mathop{\mathrm{Iso}}\nolimits}
\def\Ker{\mathop{\mathrm{Ker}}\nolimits}
\def\Lip{\mathop{\mathrm{Lip}}\nolimits}
\def\rank{\mathop{\mathrm{rank}}\limits}
\def\Ran{\mathop{\mathrm{Ran}}\nolimits}
\def\Re{\mathop{\mathrm{Re}}\nolimits}
\def\Res{\mathop{\mathrm{Res}}\nolimits}
\def\res{\mathop{\mathrm{res}}\limits}
\def\sign{\mathop{\mathrm{sign}}\nolimits}
\def\span{\mathop{\mathrm{span}}\nolimits}
\def\supp{\mathop{\mathrm{supp}}\nolimits}
\def\Tr{\mathop{\mathrm{Tr}}\nolimits}
\def\BBox{\hspace{1mm}\vrule height6pt width5.5pt depth0pt \hspace{6pt}}
\def\where{\mathop{\mathrm{where}}\nolimits}
\def\as{\mathop{\mathrm{as}}\nolimits}


\newcommand\nh[2]{\widehat{#1}\vphantom{#1}^{(#2)}}
\def\dia{\diamond}

\def\Oplus{\bigoplus\nolimits}



\def\qqq{\qquad}
\def\qq{\quad}
\let\ge\geqslant
\let\le\leqslant
\let\geq\geqslant
\let\leq\leqslant
\newcommand{\ca}{\begin{cases}}
\newcommand{\ac}{\end{cases}}
\newcommand{\ma}{\begin{pmatrix}}
\newcommand{\am}{\end{pmatrix}}
\renewcommand{\[}{\begin{equation}}
\renewcommand{\]}{\end{equation}}
\def\eq{\begin{equation}}
\def\qe{\end{equation}}
\def\[{\begin{equation}}
\def\bu{\bullet}

\title[{Resonances  for 1d Stark operators}]
{Resonances  for 1d Stark  operators}

\date{\today}
\author[Evgeny Korotyaev]{Evgeny L. Korotyaev}
\address{Saint-Petersburg State University, Universitetskaya nab.
7/9, St. Petersburg, 199034, Russia, \ korotyaev@gmail.com, \
e.korotyaev@spbu.ru}


\begin{abstract}
\no We consider the Stark operator perturbed by a compactly
supported potential (of a certain class)  on the real line. We
prove the following results:
(a) upper and lower bounds on the number of resonances in
complex discs with large radii,
  (b) the trace formula  in
terms of resonances only, (c)  all resonances determine the potential uniquely.
\end{abstract}

\maketitle

{\it Dedicated to the memory of Professor Viktor Havin,
 (St.Petersburg, 1933-2015)}

\section {Introduction and main results}
\setcounter{equation}{0}

\subsection{Introduction}
We consider the operator $H=H_0+V$ acting on
$L^2(\R)$, where the  unperturbed  operator $H_0$ is   the  Stark
operator    given by
\[
\label{a.1} H_0=-{d^2\/dx^2}+  x.
\]
Here  $x$ is an external electric field and the potential $V=V(x),
x\in \R$ is real and satisfies

\medskip

\no {\bf Condition V}. {\it The potential  $V\in L_{real}^2(\R)$
and $\supp V\ss [0,\g]$ for some $\g>0$.}

\medskip

Our main results devote  to the asymptotics of the number of
resonances in large discs
and an inverse problem in terms of resonances (all resonances
determine the potential uniquely).
Under Condition V the operator $V(H_0-i)^{-1}$ is compact (see Lemma
\ref{TP}). Then the operators $H_0$ and $H$ are self-adjoint on the
same domain  and $C_0^\infty(\R)$ is a core for both $H_0$ and $H$.
The spectrum of both $H_0$ and $H$ is purely absolutely continuous
and covers the real line $\R$ (see Avron-Herbst \cite{AH77} and
Herbst \cite{H77}).

The Stark effect is the shifting and splitting of spectral lines of atoms and molecules due to presence of an external electric field.
The effect is named after Stark, who discovered it in 1913.
The Stark effect has been of marginal benefit in the analysis of atomic spectra, but has been a major tool for molecular rotational spectra.
The perturbation theory for the Stark effect has some problems.
In absence of an electric field, states of atoms and molecules are   square-integrable. In the presence of an electric field, states of atoms and molecules are not square-integrable and they becomes resonances of finite width.   For weak fields low lying states can be regarded as bound, but   for all other cases we need to calculate  resonances and the corresponding states, which are not square-integrable.

It is well known that the wave operators $W_\pm$ for
the pair $H_0, H$ given by
$$
W_\pm=s-\lim e^{itH}e^{-itH_0} \qqq as \qqq t\to \pm\iy,
$$
exist and  are unitary (even under much less restrictive assumptions
on the potential than considered here, see \cite{AH77,H77}).  Thus
the scattering operator $S=W_+^*W_-$ is unitary. The operators $H_0$
and $S$ commute and thus are simultaneously diagonalizable:
\[
\lb{DLH0}
L^2(\R)=\int_\R^\oplus \mH_\l d\l,\qqq H_0=\int_\R^\oplus\l I_\l d\l,\qqq
S=\int_\R^\oplus S(\l)d\l;
\]
here $I_\l$ is the identity in the fiber space $\mH_\l=\C$ and $S(\l)$ is the scattering matrix (which is a  scalar function in $\l\in \R$ for our case)  for the pair $H_0, H$ (see Yajima \cite{Y81}).

\subsection{Determinants}
The main objects studied in the present paper are the resonances of
$H$ and the scattering matrix $S(\l)=   e^{-2\pi i\f_{sc}(\l)},
\l\in \R$  for the pair  $H_0, H$, where $\phi_{sc}(\l)$ is the
scattering phase (or the spectral shift function in the terminology
associated with the trace formula). In order to study resonances
 we chose an
approach where a central role is played by the {\em Fredholm
determinant}. More precisely, we set
\[
\lb{Y0}
\begin{aligned}
&V=|V|^{1\/2}V^{1\/2},\qqq V^{1\/2}=|V|^{1\/2} \sign V,
\\
& R_0(\l)=(H_0-\l)^{-1},\qqq  Y_0(\l)=|V|^{1\/2}\, R_0(\l)\,
V^{1\/2}, \qqq \l \in {\C}_\pm.
\end{aligned}
\]
Here $\C_\pm=\{\l\in\C: \pm\Im \l>0\}$ denote the upper and lower
half plane and $\l$ is a spectral parameter.
 We shortly describe standard properties of the operator-valued function $Y_0$, which we will prove in Section 2.
 We shall show
 that each
operator $Y_0(\l), \Im \l \ne 0,$ is trace class and thus we can
define the determinant:
\[
\label{a.2} D_\pm(\l)=\det (I +Y_0(\l)),\qqq \l\in\C_\pm.
\]
Moreover, we show that the function $D_\pm(\l), \l\in \C_\pm$ is
analytic in $ \C_\pm$, continuous up to the real line and
$D_\pm(\l)\ne 0$ for all $\l\in \ol\C_\pm$. Furthermore,  the
function  $D_\pm$ satisfies
\[
\lb{aD}
 D_\pm(\l)=1+O(\l^{-a})\qqq   as \qqq |\l| \to \infty, \qq
 \l\in\ol\C_\pm,
\]
for any fixed  $a\in (0,{1\/2})$, uniformly with respect to $\arg \l
\in [0,\pm\pi]$.  Thus we can define the branch  $\log D_\pm(\l), \l\in \C_\pm$ by
$
\log D_\pm(\l)=O(\l^{-a})$ as
$ |\l| \to \infty, \l\in \C_\pm$.
For each $\l\in \R$ the following
identities hold true:
\[
\label{S1}
 S(\l)={\ol D_+(\l+i0)\/D_+(\l+i0)}={D_-(\l-i0)\/D_+(\l+i0)}=
 e^{-2\pi i\f_{sc}(\l)},
\]
where  $\f_{sc}={1\/\pi}\arg D_+(\l+i0)$ is the scattering phase.
Thus the standard arguments give that the function $S(\l)$, defined
by \er{DLH0}, is continuous in $\l\in \R$. The basic properties
 of determinants (see below \er{B1c}, \er{2.2}) give the identity
\[
\lb{DDx} \ol D_+(\l)= {\det (I+Y_0(\l)^*)}=\det
(I+Y_0(\ol\l))=D_-(\ol\l) \qqq \forall\qq \l\in \C_+.
\]
Due to this identity it is enough to consider $ D_+$ or $ D_-$.  Our first preliminary theorem describes the
Fredholm determinant $D_+(\l)$  and its
asymptotics at high energy.

\no {\bf Condition C}. {\it The potential  $V$  satisfies
Condition V  and the
restriction of $V$ on the interval $(0,\g)$ is absolutely
continuous.}

\medskip

\begin{theorem}
\lb{T1} Let $V$ satisfy Condition ${\rm C}$ and and let $\t<1$. Then the function $\log D_+(\l)$  is analytic in $ \C_+$, continuous up
to the real line and satisfies
\[
\lb{aD1}
\begin{aligned}
 \log D_+(\l)={iV_0\/2\sqrt{\l}}+{O(1)\/\l^\t} \qq
\qq as \ |\l|\to \iy,\ \l\in \ol\C_+,\\
\end{aligned}
\]
where $V_0=\int_{\R}V(x)dx$, uniformly with respect to $\arg
\l \in [0,\pi]$.
In particular,
\[
\begin{aligned}
\label{S3} &\f_{sc}(\l)={V_0\/2\pi \sqrt{\l}}+{O(1)\/\l^\t}\qq as
\qqq
\l \to +\infty, \\
&\f_{sc}(\l)=O(1/\l^\t)\qqq \qqq as \qqq \l \to -\infty.
\end{aligned}
\]

\end{theorem}

\medskip

{\bf Remark.} 1)
Under our  assumptions on $V$ the proof of the asymptotic
expansion \er{aD1} is a bit technical.
If $V$ is,  e.g. in the Schwartz class, then the expansion becomes
much easier and higher order terms can  be derived as well (see \cite{IS15}), similar to the 3-dim  case in \cite{KP03}.

\subsection{Resonances}
In order to describe resonances  we recall the definition of
order  and type.

{\bf Definition.} {\it The entire function $f$ is of order $\b$ if
\[
\lb{dor}
\limsup\limits_{r\to \iy} {\log \log M(r)\/\log r}=\b ,
\]
where $M(r)=\sup _{|z|=r}|f(z)|$.  The function $f$ of positive
order $\b>0$ is of type $a\ge 0$ if
\[
\lb{dtyp}
\limsup\limits_{r\to \iy} {\log M(r)\/r^\b}=a.
\]
}


Under Condition V we will obtain an analytic continuation of
$D_+(\l), \l\in \C_+$  to the entire complex plane and information
on its zeros and obtain upper bounds on the number of
resonances  of the operator $H$.  We denote by $(\l_n)_1^{\iy}$ the sequence of zeros
in $\C_-$ of $D_+$ (counting multiplicities), arranged such that
\[
\lb{Ln}
0<|\l_1|\leq |\l_2|\leq |\l_2|\leq \dots.
\]
By definition, a zero $\l_n\in \C_-$ of $D_+$ is called a resonance.
The multiplicity of the resonance is the multiplicity of
  the corresponding zero of $D_+$. In order to obtain lower bounds on the number of
resonances  we assume that  the potential $V$
satisfies Condition C with $V\ne0$.
More precisely

\begin{theorem}
\lb{T2} Let  $V$ satisfy Condition V. Then $D_\pm(\l), \l\in
\C_\pm$ has  an analytic extension into the whole complex plane
and satisfies
\begin{equation}
\label{ED}
 |D_\pm(\l)|\leq C_0e^{{4\/3}|\l|^{3\/2}} \quad \forall \ \l\in \C,
\end{equation}
for some constant $C_0$.  Furthermore, by \er{S1}, the S-matrix
$S(\l), \l\in \R$ has an analytic extension into the whole upper
half plane $\C_+$ and a meromorphic  extension into the whole lower
half plane $\C_-$. The zeros of $S(\l),\l\in \C_+$ coincide with the
zeros of $D_-$ and the poles of $S(\l),\l\in \C_-$ are precisely
the zeros of $D_+$. Let, in addition,  $V$ satisfy Condition {\rm C} and $V(0)\ne 0$. Then $D_\pm$
 is an entire function of order ${3\/2}$ and  type ${4\/3}$.

\end{theorem}

\noindent {\bf Remark.} 1) From \er{ED} we deduce that
$D_+(\l)$ has the Hadamar factorization:
\[
\label{HFD}
 D_+(\l)=D_+(0)e^{p\l}\lim_{r\to +\infty}\prod_{|\l_n|\leq
r}\left(1-{\l\/\l_n}\right) e^{{\l\/\l_n}}, \quad \l\in  \C,
\]
uniformly on any compact subset of ${\C}$, where the constant $p$
satisfies
\[
\label{S1mathcalPk} p={D_+'(0)\/D_+(0)},\qqq \Im p=\pi \f_{sc}'(0),
\]
with $\f_{sc}(\cdot)$ defined in (\ref{S1}).

2) By \er{HFD}, the operator $H$ has an infinite number  of
resonances.

3) Due to \er{S1} the resonances  are the zeros $\l_n\in \C_-, n\ge 1$ of $D_+$ (and the poles of $S(\l)$ with the
same multiplicity) in
$\C_-$ labeled according to   \er{Ln}. The zeros of the S-matrix
$S(\l)$ are the zeros of $D_-$ in $\C_-$ given by $\ol\l_n\in \C_-,
n\ge 1$.

\medskip

Denote by $\cN(r,f)$  the number of zeros (counted according to multiplicity) of
$f$ having modulus $\leq r$.
The Lindel\"of Theorem jointly with Theorem \ref{T2}  applied to the Fredholm  determinant of the perturbed  Stark operator gives

\begin{corollary}
\lb{T3} Let  $V$ satisfy Condition {\rm V}. Then the entire function
$D_+$ satisfies
\[
\label{N1}
  \cN(r,D_+)\le C_1 r^{3\/2}
\]
for $r>0$ sufficiently large and for some  positive constant $C_1$.

Let in addition $V$ satisfy Condition {\rm C} and $V(0)\ne 0$. Then there is a sequence of positive numbers $r_j, \, j \in
\N,$ tending to $\infty$ and a positive constant $C_0>0$ such that
\[
\label{N2}
  \cN(r_j,D_+)\ge C_0 r_j^{3\/2}  \qqq \forall \ j \in \N.
  \]
\end{corollary}

{\bf Remark.} 1) We emphasize that the Hadamard factorization of
$D_+(\l)$  in \er{HFD} crucially depends on determining its  order
and type (which are equal to the order and type of the squared  Airy
function, which gives the generalized eigenfunctions for the
unperturbed Stark operator). In structure the Hadamard factorization
 looks similar to the factorization for the Schr\"odinger operator
$-\D+V$ in $\R^3$, see e.g. \cite{Z87}, \cite{SZ91}, \cite{B01}. In
fact, this is closely connected to counting the number of
resonances,  by a result of Lindel\"of (see \cite{L03}), which is
contained in Boas's book, see p.25 in \cite{Bo54} and in  Section 5.

Thus we see that the perturbed Stark operator $H$ on the real  line
has much more resonances than the corresponding Schr\"odinger
operator on the real line. In fact, the number of resonances (i.e.,
$\cN(r,D_+)\le Cr^{3\/2}$ in the disc $\{|\l_n|\le r\}$ ) of the
perturbed Stark operator $H$  on the real line corresponds to the
one for the Schr\"odinger operator on $\R^3$.
Recall that for the Schr\"odinger operator on $\R^3$ the number
$\cN_\r$ of resonances $\l_n$ in the disc $\{|\l_n|\le \r^2\}$ has
the bound $\cN_\r\le C\r^3=Cr^{3\/2}$ at $r=\r^2$. This explains the
similarity in the Hadamard factorization.

Our next corollary concerns the trace formula in terms of resonances.
So far, trace formulas for one-dimensional  Schr\"odinger operators
in terms of resonances have only been  determined in \cite{K04}.
Here we also follow the approach in \cite{K04}.

\begin{corollary}
\lb{T4} Let  $V$ satisfy Condition {\rm V} and let $R(\l)=(H-\l)^{-1}, \Im \l\ne 0$. Then  the following identity (the trace formula) holds true:
\[
\lb{Tr3} \Tr\, \big(R_0(\l)-R(\l)\big)=p+\sum_{n\ge 1}{\l\/
\l_n(\l-\l_n)},
\]
 where  the series converges absolutely and uniformly on any
compact set of $\C\sm \{\l_n, n\ge 1\}$.
\end{corollary}

{\bf Remark.} We discuss  trace formulas in Section 5 and we will show the following identity:
\[
\lb{BW} \f_{sc}'(\l)=\f_{sc}'(0)+{\l\/\pi} \Im \sum_{n\ge
1}{1\/\l_n(\l-\l_n)}\qq \forall\ \l\in \R,
\]
uniformly on any compact subset of ${\R}$. Note that the identity \er{BW}  is a Breit-Wigner
type formula for resonances (see p. 53 of \cite{RS78}).

We discuss  now inverse resonance problems. We show that all
resonances determine the potential uniquely. It is a first result
about inverse resonance problems for perturbed Stark operators.

\begin{theorem}
\lb{T6} Let the perturbed Stark operators $H_j=H_0+V_j, j=1,2$ act on $L^2(\R)$ and let each potential $V_j, j=1,2$ satisfy Condition {\rm V}.
Assume that $H_1$ and $H_2$ have the same resonances. Then
$V_1=V_2$.

\end{theorem}

{\bf Remark.} In the case of Schr\"odinger operator with a compactly supported potential on the half-line all resonances determine the potential  \cite{K04}. In the case of the real line
all resonances do not determine the potential  \cite{K05}.

\subsection{Brief overview}
Concerning previous results on resonances, we recall that from a
physicists point of view, they were first studied by Regge
\cite{R58}. Since then, properties of  resonances have been the
object of intense study and we refer to   \cite{SZ91} for the
mathematical approach in the multi-dimensional case  and references
given there.

A lot of papers are devoted to resonances of the one-dimensional
Schr\"odinger operator, see Froese \cite{F97}, Korotyaev \cite{K04},
Simon \cite{S00}, Zworski \cite{Z87} and references given there. We
recall that Zworski \cite{Z87} obtained the first results about the
asymptotic distribution of resonances for the Schr\"odinger operator
with compactly supported potentials on the real line (this result is
sharper than Corollary \ref{T3} in the present paper).
 Inverse problems (characterization, recovering,
uniqueness) in terms of resonances were solved by Korotyaev for a
Schr\"odinger operator with a compactly supported potential on the
real line \cite{K05} and the half-line \cite{K04}, see also Zworski
\cite{Z02}, Brown-Knowles-Weikard \cite{BKW03} concerning the
uniqueness.

Next, we mention some results  for  one-dimensional perturbed Stark
operators. The  one-dimensional scattering theory was considered by
Rejto-Sinha \cite{RS76}, Jensen \cite{J85},   Liu  \cite{L93}. The
one-dimensional inverse scattering problem is  studied by
Calogero-Degasperis \cite{CD78}, Graffi-Harrell \cite{GH82},
Kachalov-Kurylev \cite{KK91}, Kristensson \cite{K86}, Lin-Qian-Zhang
\cite{LQZ89}. There are a lot of results about the resonances of the
one-dimensional perturbed Stark operator, where the dilation
analyticity techniques are used, see e.g., \cite{H79}, \cite{J89}
and \cite{CFKS} and references therein. Note that compactly
supported potentials are not treated in these papers.

We  mention also interesting results about resonances for
one-dimensional Stark-Wannier operators $-{d^2\/dx^2}+\ve x+V_\pi$,
where the constant $\ve>0$ is the electric field strength and $V_\pi$ is
the real periodic potential: Agler-Froese \cite{AF85},
Grecchi-Sacchetti \cite{GS97}, Herbst-Howland \cite{HH79},  Jensen
\cite{J86}. In the case $\ve=0$, the resonances for  one-dimensional
operators $-{d^2\/dx^2}+V_\pi+V$, where $V$ is a compactly supported
potential were considered  by Firsova \cite{F84}, Korotyaev
\cite{K11}, Korotyaev-Schmidt \cite{KS12}.

\subsection{Plan of the paper}
In Section 2 we recall well known results on the spectral
representation of the Stark operator in a form useful for our
approach and obtain basic estimates on $Y_0(\l)$,  using Privalov's
Lemma. Section 3 contains the stationary representation of the
scattering matrix and the proof of Theorem \ref{T1}. Section 4
establishes the analytic continuation of $D_\pm(\l),\l\in \C_\pm$
and the meromorphic continuation of $S(\l)$ and gives the crucial
estimates on order and type leading to Corollary \ref{T3}. In
Section 5 we prove Theorem \ref{T2} and Theorem \ref{T4}.
 The
Appendix contains technical estimates needed in the proof of Lemma
2.4 which is crucial to obtain the asymptotic expansion \eqref{aD1}.

\medskip

\section {Unperturbed Stark operators}
\setcounter{equation}{0}

\medskip

\subsection{The well-known facts.}
 We denote by $C$
various possibly different constants whose values are immaterial in
our constructions. By $\cB$ and $\cB_\iy$ we denote the classes of
bounded and compact operators, respectively. Let $\cB_1$ and $\cB_2$
be the trace and the Hilbert-Schmidt class equipped with the norm
$\|\cdot \|_{\cB_1}$ and $ \|\cdot \|_{\cB_2}$, respectively. We
recall some well known facts. Let $A, B\in \cB$ and $AB,
BA\in\cB_1$. Then
\[
\lb{2.1} \Tr AB=\Tr BA,
\]
\[
\lb{2.2} \det (I+ AB)=\det (I+BA),
\]
\[
\lb{B10}
{\rm the\ mapping}\qq X\to \det (I+ X) \qq {\rm is\  continuous\ on}\ \cB_1,
\]
\[
\lb{B1}
|\det (I+ X)|\le e^{\|X\|_{\cB_1}},
\]
\[
\lb{B1c}
\ol{\det (I+ X)}=\det (I+ X^*),
\]
\[
\lb{B12}
|\det (I+ X)-\det (I+Y)|\le \|X-Y\|_{\cB_1}e^{1+\|X\|_{\cB_1}+\|Y\|_{\cB_1}},
\]
for all $X,Y\in \cB_1$, see e.g., Sect. 3. in the book \cite{S05}.
Let  the operator-valued function $\O :\cD\to \cB_1$ be analytic for
some domain $\cD\ss\C$ and $(I+\O (z))^{-1}\in \cB$ for any $z\in
\cD$. Then the function $F(z)=\det (I+\O (z))$ satisfies
\[
\lb{2.3} F'(z)= F(z)\Tr (I+\O (z))^{-1}\O '(z),\ \   \ \ \ z\in \cD.
\]
Recall that the kernels of the operators $(-\D-\l)^{-1}$ and $e^{it
\D}$ on $L^2(\R)$ have the form
\[
\lb{R0} (-\D-\l)^{-1}(x,y)={i\/2\sqrt \l}e^{i\sqrt \l|x-y|},\qqq \ \
\l\in \C\sm \ol\R_+, \qq \sqrt \l\in \C_+,
\]
\[
\lb{eit} (e^{it\D})(x,y)={e^{-i\pi/4}\/ \sqrt{4\pi
t}}e^{i|x-y|^2/4t}, \ \ \ \ \ \ \ \ \ t\ne 0,
\]
$x,y\in \R^1$. We need the identities for the Stark operator
$H_0=-{d^2\/dx^2}+x$ from \cite{AH77} given by
\[
\lb{E2} e^{-itH_0}=e^{-itx}e^{it\D}e^{-i\pa t^2}e^{-i{t^3\/3}},\qqq
 \qq  \forall\  t\in \R.
\]
where $\pa=-i{d \/d x}, \D={d^2 \/d x^2}$. The free resolvent
$R_0(\l)=(H_0-\l)^{-1}$ and its kernel $R_0(x,y,\l)$ satisfy
\[
\lb{E3} \begin{aligned}
 R_0(\l)=i\int_0^\iy  e^{-it(H_0-\l)}dt=i\int_0^\iy
e^{-itx}e^{-i\pa t^2}e^{it\D}e^{it\l-i{t^3\/3}}dt,\\
R_0(x,y,\l)={e^{i{\pi\/4}}\/\sqrt{4\pi}}\int_0^\iy
e^{-itx}e^{{i\/4t}|x-t^2-y|^2}e^{-i{t^3\/3}+it\l}{dt\/t^{1/2}},
\end{aligned}
\]
for $x,y\in \R$ and $\l\in \C_+$. We introduce the resolvent $R(\l)$
for $H$ and operators  $Y, J$ by
$$
R(\l)=(H-\l)^{-1}, \qq Y(\l) = |V|^{1\/2}R(\l)V^{1\/2},
\qq J(\l)=I-Y(\l), \qq  J_0(\l)=I+Y_0(\l),
$$
for $\l\in \C_\pm$ and recall that $Y_0(\l) =
|V|^{1\/2}R_0(\l)V^{1\/2}$. Below we will use the identities
\[
\lb{2.5}
\begin{aligned}
 \qqq R=R_0-R_0VR,\ \ \ J(\l)J_0(\l)=I,\ \ \ \ \ \l\in \C_\pm.
\end{aligned}
\]

\subsection{The spectral representation for $H_0$.}
We will need
some facts concerning the spectral decomposition of the Stark
operator $H_0$, which we denote by
$$
E_0(\l)=\c(\l-H_0),\qq \l\in \R,\qqq  \where \qqq \c(\l)=\ca 1 &
\l\ge 0\\ 0 &\l<0\ac.
$$
Now we recall formulae for $E_0(\l)$ due to \cite{AH77}. Let $\vp$
be  the multiplication operator by the function $\vp(k)=e^{ik^3/3},
k\in \R$. Then
\[
\lb{E1} H_0=-{d^2\/dx^2}+x=\vp(\pa)^* (\ x \ )\vp(\pa),\qqq \pa=-i{d\/dx}.
\]
Let $U: f\mapsto \wt f$ be the unitary transformation on $L^2(\R)$,
which can be defined on $L^1(\R)\cap L^2(\R)$ by the explicit
formula
\[
\lb{UT} \wt
f(p)=(Uf)(p)={1\/\sqrt \pi}\int_{\R} {\rm Ai}(x-p)f(x)dx,
\]
where ${\rm Ai}(\cdot)$ is the Airy function:
\[
\label{Ai01} {\rm  Ai}(z)={1\/\pi}\int_{0}^\iy \cos
\rt({t^3\/3}+tz\rt)dt,\qqq \forall \ \ z\in \R.
\]
The unitary transformation \er{UT} carries $H_0$  over  into
multiplication by $p$ in $L^2(\R,dp)$:
\[
\label{c.3} (UH_0U^* \wt f)(p) =p\wt f(p),\qqq \wt f\in \mD(p).
\]
Thus, for any $f\in L^2(\R)$, the quadratic form of $E_0(\l)$ can be presented as
\[
\label{c.4} \lan E_0(\l)f,f\ran =\int_{p<\l}|\wt f(p)|^2 dp \qqq \forall \
\l\in \R,
\]
where $\lan \cdot, \cdot\ran $ is the scalar product in $L^2(\R)$.
Differentiation with respect to $\l$ gives
$$
\frac{d}{d\l}\lan E_0(\l)f,f\ran =|G(\l)f|^2, \qq f\in L^1(\R)\cap L^2(\R),
$$
where $G(\l):L^1(\R)\cap L^2(\R)\to  \C$ is given by
\begin{equation}
\label{defG} G(\l)f:=U(\l)f={1\/\sqrt \pi}\int_{\R} {\rm
Ai}(x-\l)f(x)dx \qqq \forall \ \l\in \R.
\end{equation}

The Airy function Ai$(z), z\in \C$ is entire and satisfies (see (4.01)-(4.05)
from \cite{O74}):
\[
\label{Ai1} {\rm Ai}''(z)=z {\rm Ai}(z),
\]
\[
\label{Ai02} {\rm  Ai}(-z)=e^{i\pi/3}{\rm
Ai}(ze^{i\pi/3})+e^{-i\pi/3}{\rm  Ai}(ze^{-i\pi/3}),
\]
and it obeys the following asymptotics, as $|z|\to \iy$ uniformly in
$\arg z$ for any fixed $\ve>0$:

\[
\label{Ai2}
\begin{aligned}
{\rm Ai}(z)={1\/2}z^{-{1\/4}}e^{-{2\/3}z^{3\/2}}\rt(1+O(z^{-{3\/2}})\rt), \qq &
{\text if} \qq |\arg z|<\pi-\ve,\\
{\rm Ai}(-z)=z^{-{1\/4}}\rt[ \sin \vt+O(z^{-{3\/2}}e^{|\Im
\vt|})\rt], \qq & {\text if} \qq  \ |\arg z|\le \ve, \qq
\vt={2\/3}z^{3\/2}+{\pi \/4}.
\end{aligned}
\]

Introduce the space $L^p(\R)$ equipped by the norm
$\|f\|_p=(\int_\R|f(x)|^p dx)^{1\/p}\ge 0$. For
the scalar product  $\lan f,f\ran $ in $L^2(\R)$ and we have  $\lan f,f\ran
=\|f\|^2=\|f\|_2^2$. Define the linear functional $\P(\l):
L^2(\R)\to \C,$ by
\[
\lb{Pf} \P(\l)f=\lan f, G(\l)|V|^{1\/2}\ran ={1\/\sqrt \pi}\int_{\R} {\rm
Ai}(x-\l)|V(x)|^{1\/2}f(x)dx \qq \forall \ \l\in
\R,
\]
which is bounded, since ${\rm Ai}(\cdot)\in L^\iy(\R)$ by \er{Ai2}
and $|V|^{1\/2}\in L^2(\R)$ under Condition V.

\begin{lemma}
\label{TP} Let the potential $V$ satisfy Condition {\rm V} and let $
\t\in [0,1], \l,\m\in \R, \ |\l-\m|\le 1$. Then

i) The linear functional $ \P(\l): L^2(\R)\to \C, $ defined by
\er{Pf} is bounded and  satisfies
\[
\label{PV1} \|\P(\l)\|\le {C\/(1+|\l|)^{1\/4}},
\]
\[
\label{PV1g} \|\P(\l)-\P(\m)\|\le
{C|\l-\m|^\t\/(1+|\l|)^{{1\/4}-{\t\/2}}}.
\]
 ii) Moreover, $\O(\l)=\P(\l)^*\P(\l), \l\in \R$ is a rank
one operator on $L^2(\R)$ satisfying
\[
\label{O1} \|\O(\l)\|_{\cB_1}\le {C\/(1+|\l|)^{1\/2}},
\]
\[
\label{O2} \|\O(\l)-\O(\m)\|_{\cB_1}\le
{C|\l-\m|^\t\/(1+|\l|)^{{1-\t\/2}}}.
\]
Here the constant $C$ in \er{PV1}-\er{O2}  depends on $V$ only.

iii) Let $q, q_1$ be multiplication operators by functions $q,
q_1\in L^2(\R)$ respectively. Then for all $\Im z\ne 0$ the
following holds true:
\[
\label{e.1} qR_0(z), \qq qR(z)\in \cB_2,\qq {\and }\qqq
qR_0(z)q_1,\qq qR(z)q_1 \in \cB_1,
\]
and in particular, $VR_0(z) \in \cB_2$.

\end{lemma}
{\bf Proof.} In order to prove the lemma we need a following simple
estimate:
\[
\lb{w1}
\begin{aligned}
\int_{\R}{|V(x)|dx\/(1+|x-\l|)^{a}}\le \|V\|_2{\g(1+\g)^a\/(1+{|\l|\/2})^{a}},
\qqq \forall \ (a,\l)\in [0,{1\/2}]\ts \C.
\end{aligned}
\]
We have
\[
\int_{\R}{|V(x)|dx\/(1+|x-\l|)^{a}}\le \|V\|_2 J(\l),\qqq
J(\l)^2=\int_0^\g{dx\/(1+|x-\l|)^{2a}}.
\]
We will estimate $J(\l)$:
firstly, if $|\l|\le 2\g$, then $J^2(\l)\le \g$;
secondly, if $|\l|\ge 2\g$, then
$$
J(\l)^2\le \int_0^\g{dx\/(1+{|\l|\/2})^{2a}}={\g\/(1+{|\l|\/2})^{2a}}\le {2\g \/(1+|\l|)^{2a}},
$$
which gives \er{w1}.

i)  Let  $f\in L^2(\R)$. The asymptotics \er{Ai2} and
the estimate \er{w1} imply
\[
\begin{aligned}
|\P(\l)f|^2={1\/\pi}\rt|\int_{\R}{\rm
Ai}(x-\l)|V(x)|^{1\/2}f(x)dx\rt|^2 \le
C\rt(\int_{\R}{|V(x)|^{1\/2}|f(x)|dx\/(1+|x-\l|)^{1\/4}}\rt)^2\\
\le
C \|f\|^2\int_{\R}{|V(x)|dx\/(1+|x-\l|)^{1\/2}} \le {C_1
\|f\|^2\/(1+|\l|)^{1\/2}},
\end{aligned}
\]
for some constants $C, C_1$. This yields \er{PV1}. Next, we show \er{PV1g}. Asymptotics \er{Ai2}
entails
\[
|{\rm Ai}(\l)-{\rm Ai}(\m)|\le {C|\l-\m|^\t\/(1+|\l|)^{{1\/4}-{\t\/2}}},\qqq
 \forall \ \l,\m\in \R, \ |\l-\m|\le 1.
\]
Then similar arguments as above give (with $t={1\/4}-{\t\/2}$)
\[
\begin{aligned}
|\P(\l)f-\P(\m)f|&={1\/\sqrt \pi}\rt|\int_{\R} \big({\rm
Ai}(x-\l)-{\rm Ai}(x-\m)\big)|V(x)|^{1\/2}f(x)dx\rt|
\\
&\le
C_0|\l-\m|^\t\int_{\R}{|V(x)|^{1\/2}|f(x)|dx\/(1+|x-\l|)^{{1\/4}-{\t\/2}}}
        \le {C_2|\l-\m|^\t \/(1+|\l|)^{{1\/4}-{\t\/2}}}\|f\|,
\end{aligned}
\]
for some constants $C, C_1$. This yields \er{PV1g}. The results of ii) follow from i).

iii) For $\Im z\ne 0$ and  $q\in L^2(\R)$ due to \er{Ai2} we obtain
\[
\lb{R1} \||R_0(z)|^{1\/2}q\|_{\cB_2}^2=\int_{\R^2}{{\rm Ai}^2
(x-p)|q(x)|^2\/\pi |p-z|}dxdp<
\int_{\R^2}{C|q(x)|^2dxdp\/|p-z|(1+|p-x|)^{1\/2}}<\iy,
\]
which yields $qR_0(z), VR_0(z)\in \cB_2$ and the equality of the
domains $\mD(H_0)=\mD(H)$. From \er{R1} we deduce that
$$
qR_0(z)q_1=\rt(q|R_0(z)|^{1\/2} \rt)\rt(R_0(z)^{1\/2}q_1 \rt)\in
\cB_1.
$$
Combining these results with the identity $R=R_0-R_0V R$ we arrive
at \er{e.1}.
 \hfill  \BBox

\subsection{
 Estimates on $Y_0$.} In order to  estimate the
operator-valued functions $Y(\l), Y_0(\l), \l\in \ol \C_\pm$ in
terms of the trace class norm ${\cB_1}$ we need some additional
definitions.

Let $\cH$ be a Banach space equipped with the norm
$\|\cdot\|_{\cH}$. For any $\vt>0,  \t\in (0,1)$  we introduce the
Banach space  $\gX_{\vt,\t}(\R)=\gX_{\vt,\t}(\R, \cH)$ of the
functions $f: \R\to \cH$     equipped with the norm:
\[
\begin{aligned}
\|f\|_{\gX_{\vt,\t}(\R)}=\sup_{t\in \R, \ |h|\le 1}(1+|t|)^{\vt}
\rt(\|f(t)\|_{\cH}+   {\|f(t+h)-f(t)\|_{\cH}\/|h|^\t}\rt)<\iy,
\end{aligned}
\]
and the Banach space $\gX_{\vt,\t}(\C_\pm)=\gX_{\vt,\t}(\C_\pm,
\cH)$ of the functions $F: \C_\pm\to \cH$ equipped with the norm:
\[
\begin{aligned}
\lb{YY1} \|F\|_{\gX_{\vt,\t}(\C_\pm)}= \sup_{\l, \m\in \C_+,
|\l-\m|<1}(1+|\l|)^{\vt}\rt(\|F(\l)\|_{\cH}+
{\|F(\l)-F(\m)\|_{\cH}\/|\l-\m|^\t}\rt)<\iy.
\end{aligned}
\]

We recall Privalov's Lemma. Privalov actually proved his lemma for a
certain contour and  for scalar functions. Faddeev (see Lemma 3.1 in
\cite{F66}) proved a version for Hilbert space valued functions,
where the contour is the real line, see also \cite{A75} about the
Hilbert transformation.

\begin{lemma} {\bf (Privalov)}
\lb{TPr} Assume that $f: \R\to \cH$ belongs to the Banach space
$\gX_{\vt,\t}(\R, \cH)$ for some $\vt>0, \t\in (0,1)$ and for some
Banach space $\cH$. Then the function $F$ given by
\[
F(z)=\int_\R{f(t)\/t-z}dt,\qqq z\in \C_+,
\]
is analytic in $\C_+$ and continuous up to the real line. Moreover,
it is bounded as a map $f\to F$ from $\gX_{\vt,\t}(\R, \cH)$ into
$\gX_{\vt_1,\t}(\C_\pm, \cH)$, for any $ \vt_1<\vt$, and satisfies
\[
\begin{aligned}
\|F\|_{\gX_{\vt_1,\t}(\C_\pm ,\cH)}\le
C\|f\|_{\gX_{\vt,\t}(\R,\cH)},\qqq \forall \qq \vt_1<\vt,
\end{aligned}
\]
 where the constant $C=C(\vt_1,\vt,\t)$ depends on $\vt_1,\vt$ and $\t$
 only.
 \end{lemma}

We apply Privalov's Lemma \ref{TPr} to study the sandwiched
resolvents $Y_0(\cdot)$ and $Y(\cdot)$.

\begin{lemma}
\lb{T2.1}
 Let the potential $V$ satisfy Condition {\rm V} and let $\vt<{1-\t\/2},
 \ 0<\t<{1}$. Then

\no i) The operator-valued functions $I+Y_0$ and $(I+Y_0)^{-1}$ are
uniformly H\"older on $\ol\C_+$ and
\[
\label{e.2} (I+Y_0(\l))^{-1}\in \cB,\qqq \forall \ \l\in \ol \C_\pm,
\]
Moreover, let  $F=Y_0$ or $F=Y$. Then the operator-valued function
$F$  is analytic on $\C_\pm$ and continuous up to the real line in
the $\cB_1$--norm and satisfies
\[
\begin{aligned}
\lb{GSY} \|F\|_{\gX_{\vt,\t}(\C_\pm ,\cB_1)}<\iy.
\end{aligned}
\]
ii) Moreover, the function $D_\pm$ is  analytic in  $\C_\pm$ and
satisfies
\[
\lb{e.4} D_\pm-1\in \gX_{\vt, \t}(\C_\pm).
\]

\end{lemma}

{\bf Proof.} i) Introduce the Hilbert space $\gS_a, a\ge 0$ of
 Hilbert-Schmidt operators $X$ equipped with the norm
 $$
 \|X\|_{\gS_a}^2=\sum_{n\ge 1}n^{2a}\s_n^2,
 $$
where the non-negative numbers $\s_1\ge \s_2\ge \s_3...$ are
eigenvalues of the operator $(X^*X)^{1\/2}\ge 0$. Note that
$\|X\|_{\cB_1}\le C_a\|X\|_{\gS_a}$ for all $a>{1\/2}$, where the
constant is given by $C_a^2=\sum_{n\ge 1}n^{-2a}$.

In the free case the operator $Y_0(\l)$ has the form
\[
\begin{aligned}
Y_0(\l)=\int_\R{\O(t)\/t-\l}dt,\qqq \l\in \C_+,\\
\O(t)=\P(t)^*\P(t)V_S,\qqq V_S=\sign V.
\end{aligned}
\]
Here, due to Lemma \ref{TP}, the function $\O(t), t\in \R$ is a rank
one operator-valued function and, in particular,  $\O(t)\in \gS_1,
t\in \R$. Thus Lemma \ref{TP} shows that $\O\in \gX_{\vt,\t}(\R,
\gS_1)$ for any $\vt={1-\t\/2}, \ \ 0<\t<{1}$. Then Privalov's Lemma
\ref{TPr} yields \er{GSY} for $Y_0$.

The operator $I+Y_0(\l)$ is invertible for all $\l\in \C_\pm$, since
$H_0$ is self-adjoint and satisfies \er{e.1}. Moreover, it is
standard fact that the operator $I+Y_0(\l\pm i0)$ is invertible for
all $\l\in \R$. In fact, this follows from \er{PV1}, \er{PV1g}.
These remarks  and the properties of $Y_0$ imply that the
operator-valued functions $I+Y_0$ and $(I+Y_0)^{-1}$ are  analytic
in $\C_+$, continuous up to the real line and uniformly H\"older
continuous  in $\ol\C_+$.

Consider the case $F=Y$. Using the identity $(I+Y_0)Y=Y_0$ and the
properties of $Y_0$ and $I+Y_0$ described above we obtain the proof
of i).

ii) Due to \er{GSY} the function $D_\pm(\l)=\det (I+Y_0(\l))$ is
well-defined and analytic in $\C_\pm$. Moreover, using \er{B12} we
have
$$
|\det (I+ Y_0(\l))-1|\le \|Y_0(\l)\|_{\cB_1}e^{1+\|Y_0(\l)\|_{\cB_1}}
$$
and adding \er{GSY}, we obtain \er{e.4}.
\BBox

\begin{lemma} \lb{T2.4}
Let $V$ satisfy Condition {\rm V} and let $\t<{1\/2}$. Then
\[
\lb{Y2} |\Tr Y_0^n(\l)|\le  \|Y_0^n(\l)\|_{\cB_1}
 \le {C\/|\l|^{n\t}} \ \ \ \forall \ \ n\ge 1,
\]
where $\l\in \ol\C_+, |\l|>1$ and $ C=C(\t,V)$ depending on
$\t, V$. Let, in addition, $V$ satisfy Condition {\rm C}. Then
\[
\lb{Y1} \Tr Y_0(\l)={i\/2\sqrt{\l}}V_0+{O(1)\/\l}
\]
as $|\l|\to \iy, \l\in \ol\C_+$, where $V_0=\int_\R V(x)dx$
uniformly in $\arg \l\in [0,\pi]$.
\end{lemma}

\no {\bf Proof.} For  $\l\in \ol\C_+$, due to \er{GSY} we have
\[
\begin{aligned}
|\Tr Y_0^n(\l)|\le  \|Y_0^n(\l)\|_{\cB_1}\le \|Y_0(\l)\|_{\cB_1}^n \le
{C_\t\/|\l|^{n\t}},
\end{aligned}
\]
where the constant $C=C(\ve, V)$
does not depend on $\l$. We show \er{Y1} in Lemma \ref{TAY}.
 \BBox

\medskip

\section {Determinants and S-matrix}
\setcounter{equation}{0}

\medskip

\subsection{
The Determinants.} We discuss the determinant $D_\pm(\l),
\l\in\ol\C_\pm$, when the potential $V$ satisfies Condition {\rm V}.
In this case we have \er{e.1} and this  gives the identity
\[
\lb{D1} {D_\pm'(\l)\/D_\pm(\l)}=\Tr R_0(\l)VR(\l)=\Tr (R_0(\l)-R(\l)), \qqq
\l \in \C_\pm,
\]
which is well-known for large class of operators.
Due to \er{GSY} the operator-valued function
$Y_0(\l)$ attains value in $\cB_1$ and belongs to the class
$\gX_{\vt, \t}(\C_\pm, \cB_1)$ for any  $\vt<{1-\t\/2}$ and  $0<\t<1$.
Recall that we define $\log D_\pm(\l)$, by
$\log D_\pm(\l) = o(1)$ as $|\l| \to \infty, \ \l\in \C_\pm$, since
$D_\pm(\l)\ne 0$
for all $\l\in\ol\C_\pm$ and $\|Y_0(\l)\|_{\cB_1}=o(1)$ as $|\l| \to
\infty$ and there exists $r_0>0$ such that
\[
\lb{D2.12a} \sup_{\l\in \ol \C_+, |\l|\ge r_0}\|Y_0(\l)\|_{\cB_1}<{1\/2}.
\]
Then using \er{B10}-\er{B12}, \er{GSY} we obtain
\[
\lb{D2}
 \log D_\pm\in \gX_{\vt, \t}(\C_\pm), \qqq
 \forall \ \vt<{1-\t\/2}, \ 0<\t<1.
\]
It is well-known (see
\cite{RS78}) that, under the condition \er{D2.12a}, the function
$\log D_\pm(\l)$ satisfies
\[
\lb{2.12} -\log D_\pm(\l)=\sum _{n=1}^\iy {1\/n}\Tr (-Y_0(\l))^n, \ \ \
\]
for any $\l\in \ol\C_\pm, \  |\l|>r_0$, where the series converges absolutely and uniformly.
Then  using \er{2.12} and \er{D2.12a} for some $r_0>0$ and any $\t<{1\/2}$
we obtain
\[
\lb{2.13}
 |\log D_\pm(\l)+\sum _{n=1}^{N}{1\/n}\Tr (-Y_0(\l))^n|\le
{\|Y_0(\l)\|_{\cB_1}^{N+1}\/N+1} \le
{C\/|\l|^{(N+1)\t} },\ \ \ \forall  \ N\ge 0,
\]

\subsection{
The scattering matrix.} Recall that the S-matrix
$S(\l)$ is a scalar function of $\l\in \R$, acting as multiplication
in the fiber spaces $\C=\mH_\l$. Thus $|S(\l)|=1$ for all $\l\in \R$ we have
\[
 \label{Sx} S(\l)=e^{-2\pi i\f_{sc}(\l)},\qqq \l\in\R.
\]
The stationary representation for
the scattering matrix has the form  (see e.g. \cite{Y81}):
\[
\label{S}
\begin{aligned}
&S(\l)=I-2\pi i  \cA(\l), \qqq \l\in\R;\qqq \cA=\cA_0-\cA_1, \\
&\cA_0(\l)=\P(\l)V_S\P^*(\l),\qqq  \cA_1(\l)
=\P(\l)V_SY(\l+i0)\P^*(\l), \\ &\P(\l)=G(\l)|V|^{1\/2},\qq V_S=\sign
V,
\end{aligned}
\]
where $G$ is given by \er{defG}.
Note that due to \er{GSY} the operator $Y(\l\pm i0)$ is continuous
in  $\l\in\R$.
We shall represent $S(\l)$ in terms of $D_\pm(\l)$.

\begin{lemma}
\label{TA}
 Let $V$ satisfy Condition {\rm V}.  Then the scattering amplitude $\cA(\l)$ is
a continuous  scalar function of $\l\in \R$ and satisfies
\[
\label{A0i} \cA_0(\l)= \int_\R {\rm Ai}(x-\l)^2 V(x)dx\qqq \forall \
\l\in \R,
\]
\[
\label{A0} \cA_0(\cdot )\in \gX_{{1-\t\/2}, \t}(\R,\C),
\]
\[
\label{A1x} \cA_1(\cdot )\in \gX_{\ve, \t}(\R,\C)\qq \forall
\ \ve<1-\t,
\]
for any  $\t\in (0,{1})$.  Moreover, the functions $S(\l),
\f_{sc}(\l)$  are continuous in $\l\in \R$ and satisfy asymptotics
\[
\label{A2} S(\l)-1=O(\l^{-{1\/2}}),\qqq    \f_{sc}(\l)=O(\l^{-{1\/2}}) \qqq as \qq \l\to \pm \iy,
\]
\[
\label{fsc1}
\f_{sc}(\l)={1\/\pi}\arg D_+(\l),\qqq \l\in \ol\C_+,
\]
and the identities \er{S1}  which uniquely defines $\f_{sc}$ by
\er{Sx}, continuity and the asymptotics \er{A2}.
\end{lemma}

 {\bf Proof.} The definitions of $ \cA_0$ and $ \P$ (see \er{S})
give \er{A0i}. Relation \er{A0} follows from Lemma \ref{TP}.
Relation \er{A1x} follows from Lemma \ref{TP} and \er{GSY}, since
$Y(\cdot)\in \gX_{\vt, \t}(\C_\pm, \cB_1)$,  for any $\vt<{1-\t\/2},
\t\in (0,{1})$ due to \er{GSY}.

Next, we show \er{S1}. Recall that $S(\l)$ satisfies  \er{S} and
that we have the standard identity
\[
\label{Sxy}
\begin{aligned}
Y_0(\l+i0)-Y_0(\l-i0)=2\pi i\P(\l)^*\P(\l)V_S, \qq \l\in \R.
\end{aligned}
\]
Then \er{S}, \er{2.2} give
$$
\begin{aligned}
&\det S(\l)=\det \rt(I-2\pi i \P(\l)V_SJ(\l+i0)\P(\l)^*\rt)\\
&= \det \rt(I-2\pi i \P(\l)^*\P(\l)V_SJ(\l+i0)\rt) = \det
J(\l+i0)\det
(J_0(\l+i0)-2\pi i \P(\l)^*\P(\l)V_S)\\
&=\det J(\l+i0)\det
\rt(J_0(\l+i0)-Y_0(\l+i0)+Y_0(\l-i0)\rt)\\
&=\det J(\l+i0) \det J_0(\l-i0),
\end{aligned}
$$
which together with \er{2.5} yields \er{S1} since
$$
\det S(\l)=\det J(\l+i0) \det J_0(\l-i0)={\det J_0(\l-i0)\/\det
J_0(\l+i0)}={D_-(\l-i0)\/D_+(\l+i0)}={\ol D_+(\l+i0)\/D_+(\l+i0)}.
$$
 This yields \er{fsc1} and adding  the relation \er{D2} we obtain
 $\f_{sc}(\l)=O(|\l|^{-a})$ as $\l\to \pm\iy $ for any $a<{1\/2}$.
Substituting estimates  \er{GSY} and \er{PV1} into  \er{S} we obtain
$S(\l)-1=O(\l^{-{1\/2}})$ as $ \l\to \pm \iy$. As both $\f_{sc}(\l)$
and $S(\l)$ are continuous in $\l$, formula \er{Sx} determines
$\f_{sc}(\l)$ by $\f_{sc}(\l)={i\/2\pi }\log S(\l)$ and the
asymptotics $\f_{sc}(\l)=O(\l^{-{1\/2}})$ as $ \l\to \pm \iy$. All
together this proves Lemma \ref{TA}. \BBox

{\bf Proof of Theorem \ref{T1}.} Due to Lemma \ref{T2.1}
the function $D_\pm$ is analytic in $ \C_\pm$, continuous up to the
real line.
Asymptotics \er{2.13} and \er{Y1} yield \er{aD1}, which gives  \er{S3}.
\BBox

\begin{proposition}
\lb{Ttr} Let  $V$ satisfy Condition {\rm C}. Then  the following trace formulas hold true:
\[
\lb{D1x} \int_{\R}V(x)dx={2\/\pi}\int_\R  \Re{\log
D_+(\l+i0)\/\sqrt {\l+i0}}d\l,
\]
\[
\lb{D1xx} \lim_{r\to \iy} \int_{-r}^r  \Im{\log
D_+(\l+i0)\/\sqrt {\l+i0}}d\l=0.
\]
\end{proposition}

{\bf Proof.} Define a contour : $ \G_r=c_r\cup (-r,r) $, where
$c_r=\{|\l|=r\}\cap \C_+$ for large $r\to +\iy$. The function
$f(\l)=i {\log D_+(\l)\/\sqrt \l}$ is analytic in the
upper-half-plane and continuous up to the real line without zero.
This gives
$$
0=\int_{\G_r}f(\l)d\l=I_r+I_r^+,\qqq I_r^+=\int_{c_r}f(\l)d\l, \qqq
I_r=\int_{-r}^rf(\l)d\l.
$$
 Due to asymptotics \er{aD1} $f(\l)={iV_0+O(\l^{-\ve })\/2\l}$ as $|\l|\to \iy $ uniformly in $\arg \l\in [0,\pi]$ for some $\ve>0$, we obtain
\[
\lb{tr1}
I_r^+=\int_{c_r}f(\l)d\l={iV_0\/2}\int_{c_r}{d\l\/\l}+o(1)=-{\pi \/2}V_0+o(1)\qqq
\]
and
\[
\lb{tr2}
I_r=\int_{-r}^rf(\l)d\l={\pi \/2}V_0+o(1)\qqq
\]
as $r\to \iy$, and here
\[
\lb{tr3}
\Re f(\l+i0)\in L^1(\R),\qqq \Im f(\l+i0)={V_0\/2\l}+O(\l^{-\ve-1 })
\qqq \as \ \pm\l\to \iy.
\]
 Combining all relations \er{tr1}-\er{tr3} we obtain \er{D1x}-\er{D1xx}.
 \BBox

{\bf Remark.} We recall that trace formulas are important to study
 non linear equations, inverse problems, spectral theory, etc. see
\cite{DK99}, \cite{FZ71},  \cite{KK97}, \cite{LW00}  and references
therein. The complete asymptotic expansion of the scattering phase
(the spectral shift function) $\f_{sc}$ at high energies and a
sequence of trace formulas for 3-dim perturbed Stark operators were
determined by Korotyaev-Pushnitski \cite{KP03}.

\medskip



\section {Analyticity of $\cA_0$ }
\setcounter{equation}{0}

\medskip

\subsection{Estimates on Airy functions.}
In this section we assume that the potential $V$ satisfies Condition
V.
Recall that by \er{S}, the functional $\P(\l): L^2(\R)\to \C$ and its adjoint $\P^*(\l):\C\to L^2(\R)$ for $\l\in \R$ are  given by
\[
\lb{deP} \P(\l)f=\int_{\R}{\rm Ai}(x-\l)|V(x)|^{1\/2}f(x)dx, \qq
\P^*(\l)c={\rm Ai}(x-\l)|V(x)|^{1\/2}c,\qqq \l\in \R,
\]
where $(f,c)\in L^2(\R)\ts \C$. Assuming that $V$ has compact
support, these mappings  are  bounded on the real line and have
analytic extensions from $\R$ onto the whole complex plane. We
remark that $\supp V$ being bounded from below suffices to render
$\P(\l)$ analytic. To prove our estimates $\supp V$ being compact
is, however, helpful.
 Thus  $\P^*(\ol\l)$ can be identified with the function ${\rm
Ai}(\cdot-\l)|V(\cdot)|^{1\/2}$ in $L^2(\R)$, which is analytic in $\l\in \C$.

In order to estimate $\P, \P^*$, we need the asymptotics of the Airy function
from \er{Ai2}. Furthermore we have
\[
\lb{sym1}
\begin{aligned}
\ca (x-\l)^{-{1\/4}}=(-\l)^{-{1\/4}}\big(1+{O(|\l|^{-1}}  \big),\\
(x-\l)^{{3\/2}}=(-\l)^{3\/2}+{3\/2}x(-\l)^{1\/2}+{O(|\l|^{-{1\/2}}})\ac,\qqq
|\arg \l| \geq \ve,
\end{aligned}
\]
and
\[
\lb{sym1x}
\begin{aligned}
\ca (\l-x)^{-{1\/4}}=\l^{-{1\/4}}\big(1+{O(|\l|^{-1}}  \big),\\
(\l-x)^{{3\/2}}=\l^{3\/2}-{3\/2}x \l^{1\/2}+{O(|\l|^{-{1\/2}}})\ac,\qqq
|\arg \l | \leq \ve,
\end{aligned}
\]
locally uniformly in $x \in \R$, as  $|\l|\to \iy$. We will use
these estimates in order to determine asymptotics of Airy functions.
Asymptotics \er{sym1}, \er{Ai2} with $z=x-\l$  and straightforward
calculation give  the following symptotics \er{as1}, \er{as2}:

{\it i) Let $|\arg\l| \geq \ve$ and let $\z=\sqrt {-\l},\  |\arg \z|\le {\pi-\ve\/2}$. Then as $|\l| \to \infty $   one has
\[
\label{as1}
\begin{aligned}
{\rm Ai}(x-\l)^2={1\/4\z} e^{-({4\/3} \z^3+2x \z)}
\rt(1 + {O(1)\/\z}\rt),
\end{aligned}
\]
and, in particular,}
\[
\label{as2}
|{\rm Ai}(x-\l)|^2={1\/4|\l|^{1\/2}}
e^{- \Re ({4\/3} \z^3+2x \z)}\rt(1 + {O(1)\/|\l|^{1/2}} \rt).
\]

Moreover, for  the case $|\arg \l| \leq \ve$ using \er{sym1x}, we set $X= {4\/3} (\l-x)^{3\/2}$ and
obtain after short calculation
\[
X= {4\/3} \l^{3\/2} - 2x \l^{1\/2} + O(|\l|^{-{1\/2}})=\e+O(|\l|^{-{1\/2}})\qqq
as\qqq |\l|\to \iy,\qq |\arg \l| \leq \ve.
\]
where $\e={4\/3}\z^{3}-2x\z$.
Substituting this  into  \er{Ai2} for $|\arg \l| \leq \ve$, we obtain the following asymptotics  and the estimate:

ii) {\it  Let $|\arg \l| \leq \ve, \z=\sqrt \l, |\arg \z|\le{\ve\/2}$
and let $|\l|$ be sufficiently large. Then
\[
\label{as3}
\begin{aligned}
{\rm Ai}(x-\l)^2={1+ \sin \e\/2\z}+{O(e^{|\Im \e|})\/\l},
\end{aligned}
\]
where $\e={4\/3}\z^{3}-2x\z$,
and in particular, the following
slightly weaker estimate holds true:
\[
\lb{as4} |{\rm Ai}(x-\l)|^2  \leq {1\/|\l|^{1\/2}} e^{|\Im
\e|}\rt(1+{O(1)\/|\l|^{1\/2}} \rt).
\]
Note that all estimates \er{as1}-\er{as4} are locally uniform in $x$
 on bounded intervals.}

\subsection{Estimates on the Born term $\cA_0$.}
Now we are ready to study the Born term $\cA_0$.

\begin{lemma}
\lb{TA0}
Let $V$ satisfy Condition V and let $\ve >0$. Then the Born term $\cA_0(\l), \l\in \R$ given by \er{A0i} has an analytic extension from the real line  into the whole complex plane and satisfies as $|\l|\to \iy$:

i) Let   $|\arg\l| \geq \ve$ and let $-i\z=\sqrt {-\l},\qqq |\arg \z-{\pi\/2}|\le \ve$. Then
 \[
\label{aA0}
\begin{aligned}
\cA_0(\l)={i\/4\z} e^{-i{4\/3} \z^3}\int_0^\g e^{i2x \z}V(x)
\rt(1 + {O(1)\/\z}\rt)dx,
\end{aligned}
\]
and
\[
\label{aA0e}
\begin{aligned}
|\cA_0(\l)|\le {C\/|\z|}e^{{4\/3}\Re \z^3}\|V\|_{L^1(0,\g)},
\end{aligned}
\]
for some absolute constants $C$.

ii) Let   $|\arg\l| \leq \ve$ and let $\z=\sqrt \l,\  |\arg \z|\le \ve$. Then
\[
\label{aA01}
\begin{aligned}
\cA_0(\l)={1\/2\z}\int_0^\g V(x)  \rt[1+ \sin \e+{O(e^{|\Im \e|})\/\z}\rt]dx,
\end{aligned}
\]
where $\e={4\/3}\z^{3}-2x\z$    and
\[
\label{aA01x}
\begin{aligned}
|\cA_0(\l)|\le {1\/2|\z|}\int_0^\iy |V(x)|  \rt[1+ |\sin  \e|+{O(e^{|\Im \e|})\/|\z|}\rt]dx, \qq
\end{aligned}
\]
 Let, in addition,  $|\z|\ge 1+\g$ and
 $|\arg \z|\le {\pi\/6}$. Then
 \[
\label{A0g}
\begin{aligned}
|\cA_0(\l)|\le
{e^{{4\/3}|\Im \z^3|}\/2|\z|}\int_0^\g |V(x)|  \rt(1 +{O(1)\/|\z|}\rt)dx.
\end{aligned}
\]

\end{lemma}

{\bf Proof.} i) Substituting  \er{sym1}, \er{as1} in \er{A0i} we obtain
\er{aA0} and  \er{aA0e}.

ii)  Using \er{as2} we obtain \er{aA01}, which yields \er{aA01x}.
Let $\z=re^{i\f}$, where  $|\z|\ge 1+\g$ and
 $ |\f|\le {\pi\/6}$. Then we have
 $$
 |\Im \e|=|\Im  \rt({2\/3}\z^{3}-x\z\rt)|=r\rt|{2\/3}r^2\sin 3\f-x\sin \f\rt|\le {2\/3}r^3\sin 3|\f|.
 $$
Substituting the last estimate into \er{aA01} we obtain \er{A0g}.
\BBox

\section{Resonances and S-matrix}
\setcounter{equation}{0}

\subsection{Analyticity of S-matrix}
We discuss a meromorphic continuation
of the S-matrix $S(\l)$ from the real line onto the whole complex plane.

\begin{lemma}
\label{Am1} Let  $V$ satisfy Condition {\rm V} and let $\ve >0$.
Then

i) The functionals $\P(\l): L^2(\R)\to \C$ given by \er{deP},
and the mapping $\P^*(\ol\l)=\P^*(\l)$, for all $\l\in \R$
have analytic extensions from the real line into the whole complex plane and
satisfy
\[
\lb{Px1} \|\P(\l)\|^2 = \|\P^*(\ol\l)\|^2 =\int_0^\g |{\rm Ai}(x-\l)|^2
|V(x)|dx,
\]
 \[
\lb{Px1a}
\begin{aligned}
\|\P(\l)\|^2\le {C\/(1+|\l|)^{1\/2}}e^{{4\/3}|\l|^{3\/2}},
\end{aligned}
\]
for all  $\l\in \C$ and for some constant $C=C(V)$.

ii) The scattering amplitude $\cA(\l) = \cA_0(\l)-\cA_1(\l)$,
defined in \er{S} for $\l \in \R$, has an analytic extension from
the real line into the whole upper half-plane satisfying
\[
\lb{A41} |\cA_1(\l)|\le  {C_a \/ (1+|\l|)^a} \|\P(\l)\|^2 \qqq \forall
\qq \l \in \ol\C_+,
\]
for any $a\in (0,{1\/2})$ and some constant $C=C_a$ depending on
$a$.

\end{lemma}

{\bf Proof.} i) Since the Airy function is entire, the functional
$\P(\l): L^2(\R)\to \C$ given by \er{deP}, and the mapping
$\P(\l)^*$, for all $\l\in \R$ have analytic extensions
from the real line into the whole complex plane. The  identities in
\er{Px1} are obvious from the definitions of $\P, \P^*$. The proof
of \er{Px1a} is a repeatition of the proof of \er{A0g} and \er{aA0}.
In fact, substituting \er{as2}, \er{as4} into \er{Px1} we obtain
\er{Px1a}.

ii) The operators $\P, \P_1$ have analytic extensions from the real
line into the whole complex plane and operator-valued function
$Y(\l+i0)$ also has an analytic extension from the real line into
the upper-half plane. Then $\cA_1(\l)=\P(\l)V_SY(\l+i0))\P^*(\l)$
has an analytic extension from the real line into the upper-half
plane and thus $\cA$ has so. Moreover, \er{GSY}  gives
\[
    \lb{a1} |\cA_1(\l)| \leq \|\P(\l)\|^2 \| Y(\l) \| \le C
(1+|\l|)^{-a} \|\P(\l)\|^2 \qq \forall \ \l\in \ol\C_+,
\]
where $C=C(V)$ is some constant.\BBox

\begin{lemma}
\label{Am2}
Let $V$ satisfy Condition {\rm C} with $V(0)\ne0$ and let
$\l=t^2e^{i{\pi\/3}}, \ \z=e^{i{\pi\/6}}t$ as  $t\to +\iy$. Then
 \[
\label{aA0x}
\begin{aligned}
\cA_0(\l)={e^{i{\pi\/3}+{4\/3} t^3}\/4 t} \int_0^\g
e^{(-1 +i \sqrt 3)tx}V(x)
\rt(1 + {O(1)\/t}\rt)dx={e^{i{2\pi\/3}+{4\/3} t^3}\/8t^2}(V(0)+o(1)),
\end{aligned}
\]
 \[
\label{A11}
\|\P(\l)\|^2\le e^{{4\/3}t^3}{C\/t^2}\|V\|,
\]
\[
\label{A12}
\cA(\l)=
{e^{i{2\pi\/3}+{4\/3} t^3}\/8t^2}(V(0)+o(1)).
\]
\end{lemma}

{\bf Proof.}  Let $\s=1 -i \sqrt 3=2e^{-i{\pi\/3}}$ and let $ \z=e^{i{\pi\/6}}t$ as $t\to +\iy$. From  \er{aA0} we have
\[
\label{aA0xxx}
\begin{aligned}
\cA_0(\l)={ie^{-i{\pi\/6}}\/4 t} e^{{4\/3} t^3}\cI(\l),\qq
\cI(\l)&=\int_0^\g
e^{-\s tx}V(x)\rt(1 + {O(1)\/t}\rt)dx=\cI_1(\l)+\cI_2(\l),\\
\cI_1(\l)=&\int_0^{\g_1} e^{-\s tx}V(x)\rt(1 + {O(1)\/t}\rt)dx,\\
\cI_2(\l)=&\int_{\g_1}^\g e^{-\s tx}V(x)\rt(1 +
{O(1)\/t}\rt)dx=e^{-t\g_1}\|V\|_1O(1).
\end{aligned}
\]
An integration by parts yields
$$
\begin{aligned}
\int_0^\g
e^{-\s tx}V(x)dx=-{1\/\s t}e^{-\s tx}V(x)\rt|_0^\g+{1\/\s t}\int_0^\g
e^{-\s tx}V'(x)dx\\
={1\/\s t}\rt(V(0)-e^{-\s t\g }V(\g)+\int_0^\g
e^{-\s tx}V'(x)dx\rt)={V(0)-o(1)\/\s t},
\end{aligned}
$$
since since $V$ is absolutely continuous, its derivative is in
$L^1(\R)$ and for $u\in L^1(0,\iy)$ the following asymptotics  hold
true
\[
\lb{LTu} \int_0^\g e^{-tx}u(x)dx=o(1)\qqq as \qqq t\to \iy.
\]
Moreover, \er{LTu} gives $ \int_0^\g e^{-\s tx}V(x)dx=o(1)$.
Combining all these estimates we obtain \er{aA0x}.

 We consider $\|\P(\l)\|^2$ and we show \er{A11}.
Let $ \z=te^{i{\pi\/6}}$ as $t\to +\iy$. From  \er{aA0} we have
$$
\|\P(\l)\|^2=\int_0^\g |{\rm Ai}(x-\l)|^2 |V(x)|dx={1\/4 t} e^{{4\/3} t^3}\int_0^\g e^{-tx}|V(x)|(1 +o(1))dx,
$$
and
$$
\begin{aligned}
\int_0^\g e^{-tx}|V(x)|dx
&=\int_0^{\g_1}e^{-tx}|V(x)|dx+\int_{\g_1}^{\g}e^{-tx}|V(x)|dx\\
&\le \|V\|_{L^\iy(0,\g_1)}\int_0^{\g_1}
e^{-tx}dx+e^{-t\g_1}\int_{\g_1}^\g |V(x)|dx\le {1\/t}\|V\|_\iy
+e^{-t\g_1}\|V\|_1,
\end{aligned}
$$
which yields \er{A11}.  By \er{S}, the S-matrix has the form
$$
S(\l)=I-2\pi i  \cA(\l),\qqq  \cA=\cA_0-\cA_1.
$$
 Using \er{A41},\er{A11} we
obtain
 $$
 \cA(\l)=\cA_0(\l)-\cA_1(\l)={e^{i{2\pi\/3}+{4\/3} t^3}\/8t^2}\rt[(V(0)+o(1))+
 O(t^{-a})\rt],
 $$
which permits to get \er{A12}.
 \hfill  \BBox

We are now ready to prove the main  theorems.

{\bf Proof of Theorem  \ref{T2}.}
It follows  from Lemma \ref{Am1} and   \ref{TA0} that $S(\l), \l\in \R$ has an analytic extension from $\R$ into the upper half-plane ${\C_+}$
and satisfies
\[
\lb{esS}
|S(\l)|\le C_Se^{{4\/3}|\l|^{3\/2}} \qqq \forall
\ \l\in \ol\C_+,
\]
for some constant $C_S$. Invoking the identity \er{S1} we obtain $
 D_-(\l-i0) = S(\l) D_+(\l+i0),\  \l \in \R$.
Thus  it follows  from Lemma \ref{Am1}, \ref{TA0} and   \ref{T2.1} that $D_-(\l), \l\in \C_-$ has an analytic extension from $\ol{\C}_-$ into the whole complex plane ${\C}$ and satisfies
$$
\begin{aligned}
& |D_-(\l)|\le |S(\l)||D_+(\l)|\le C_Se^{{4\/3}|\l|^{3\/2}} \sup_{\l\in \C_+}|D_+(\l)|\qqq \forall
\ \l\in \ol\C_+,\\
&\sup_{\l\in \C_\pm}|D_\pm(\l)|<\iy,
\end{aligned}
$$
which together with \er{DDx}   yields \er{ED}.
Then, by \er{S1}, the S-matrix
$S(\l), \l\in \R$ has an analytic extension into the whole upper
half plane $\C_+$ and a meromorphic  extension into the whole lower
half plane $\C_-$  satisfying
\[
\lb{SDD}
S(\l)={D_-(\l)\/D_+(\l)}\qqq \forall \ \l\in \C,
\]
where the functions $D_\pm$ do not vanish in $\ol{\C}_\pm$.
The zeros of $S(\l),\l\in \C_+$ coincide with the
zeros of $D_-$ and the poles of $S(\l),\l\in \C_-$ are precisely
the zeros of $D_+$.

Let, in addition, $V$ satisfy Condition C with $V(0)\ne0$. Then
 asymptotics \er{A12} , \er{aD1}  give
 \[
\label{D11}
\begin{aligned}
D_-(t^2e^{i{\pi\/3}})={e^{i{2\pi\/3}+{4\/3} t^3}\/8
t^2}(V(0)+o(1))\qqq \as \qq \qq   t\to +\iy.
\end{aligned}
\]
Thus $D_\pm$ is an entire functions of order ${3\/2}$ and type
${4\/3}$.
\BBox

\

Now we discuss Remarks after Theorem \ref{T2}.
By Theorem  \ref{T2}, the determinants
$D_\pm(\l),\l\in \C_\pm$  extend to entire functions of order ${3\/2}$ and type
${4\/3}$.
 It is well known that in this case $D_+$ has the Hadamard
factorization \er{HFD} and \er{DD} for some $p$, see p. 22 in
\cite{Bo54}. Moreover, we have
\[
\label{DD}
 {D_+'(\l)\/D_+(\l)}=p+\sum_{n\geq1} {\l\/\l_n(\l-\l_n)},
\]
where the constant $p={D_+'(0)\/D_+(0)}$ and
the  series converges absolutely and uniformly on any
compact set of $\C\sm \{\l_n, n\ge 1\}$.
Recall that if an entire function $F$ has order $m>0$
and has zeros $z_n, n\ge 1$, then (see p. 17 in \cite{Bo54})
\[
\lb{mm1} \sum _{z_n\ne 0, n\ge 1}|z_n|^{-m_1}<\iy, \qqq \forall \
m_1>m.
\]

Using \er{HFD}, \er{DD} and differentiating $S(\l)={ \ol
D_+(\l+i0)\/D_+(\l+i0)}=e^{-2\pi i\f_{sc}(\l)}, \l\in \R$, we obtain
$$
S(\l)(-2\pi i\f_{sc}'(\l))=-S(\l)\rt({D_+'(\l+i0)\/D_+(\l+i0)}-{\ol
D_+'(\l+i0)\/\ol D_+(\l+i0)} \rt),\qqq \l\in \R,
$$
which yields
\[
2\pi i\f_{sc}'(\l)=p-\ol p+\l \sum_{n\ge 1}\rt({1\/\l_n(\l-\l_n)}-
{1\/\ol\l_n(\l-\ol\l_n)} \rt),\qq \forall\ \l\in \R,
\]
and then $2\pi i\f_{sc}'(0)=p-\ol p$. Thus we obtain \er{BW}.

Finally,  we prove  Corollary
\ref{T3} and \ref{T4}.

{\bf Proof of Corollary  \ref{T3}.}
Under the Condition V Theorem  \ref{T2} gives that  the determinant $D_+(\l),\l\in \C_+$  extends to an entire function
and satisfies \er{ED}, which yields the standard upper bound \er{N1}
(see page 16 in \cite{Bo54}).

Let $V$ satisfy the Condition C and $V(0)\ne 0$. Then
by Theorem  \ref{T2},
determinant $D_+(\l),\l\in \C_+$  is an entire function of
order ${3\/2}$ and type ${4\/3}$. Since the order of $D_+$ is not
integer, \er{N1} follows from  the Lindel\"of Theorem \cite{L03}.
We recall this theorem.

\no {\it Lindel\"of Theorem ( 1903).} {\it  Let $f$ be an entire
function of order $\b$, which is not an integer. Then \qqq i) $f$ is
of zero type iff $\cN(r,f)=o(r^\b)$,

\qqq \qq \ ii)  $f$ is of finite type iff $\cN(r,f)=O(r^\b)$. }

\medskip
\no Moreover, since the type of $D_+$ is different from
zero,  $\cN(r,D_+)$ is not $o(r^{3/2})$. This is \er{N2}. \hfill
\BBox

{\bf Proof of Corollary   \ref{T4}.} Due to \er{D1} we have
$\Tr (R(\l)-R_0(\l))={D_+'(\l)\/D_+'(\l)}$ for each $\l\in \C_+$
and adding \er{DD} we obtain    \er{Tr3}.
\BBox

\begin{proposition}
\lb{T5} Let  $V$ satisfy Condition {\rm V}. Then for any $f\in
C_0^{\infty}(\R)$ the following identity
\[
\label{Tr1} \Tr\,\left(f(H)-f(H_0)\right)= -{1\/\pi}\int_\R
f(t)\rt( \f_{sc}'(0)+{\l\/\pi} \Im \sum_{n\ge 1}{1\/\l_n(t-\l_n)}
\rt)dt,
\]
 holds true, furthermore,

\[
\lb{diff} \f_{sc}^{(m)}(0)={(m-1)!\/\pi}\Im \sum _{n\ge
1}{1\/\l_n^m},
\]
 where  the first  series converges absolutely and uniformly on any
compact set of $\C\sm \{\l_n, n\ge 1\}$.

\end{proposition}

{\bf Proof.}
Due to \er{e.1} we obtain
$R(\l)-R_0(\l)\in \cB_1$ for each $\Im \l\ne 0$. Then  the Krein
formula \cite{Kr62} for the operators $H, H_0$ and for any $f\in
C_0^{\infty}(\R )$ gives
\[
\lb{Kr1} \Tr\,\left(f(H)-f(H_0)\right)= -\int_\R
f(t)\f_{sc}'(t)dt
\]
and substituting \er{BW} into  \er{Kr1} we obtain \er{Tr1}.

Due to \er{BW} we have the identity
\[
\lb{df} \f_{sc}'(\l)=\f_{sc}'(0)+{1\/\pi} \Im \sum_{n\ge
1}\rt({1\/\l_n-\l}-{1\/\l_n} \rt),\qq \forall\ \l\in \C,
\]
uniformly on any compact subset of ${\C}\setminus \{\l_1, \l_2,
\l_3,\cdots\}$.
Differentiating $\f_{sc}'(\l)$ we then arrive at \er{diff}. \BBox

\medskip

{\bf Proof of Theorem \ref{T6}.}
There are results about inverse problems for
perturbed Stark operators $H=H_0+V$ on $L^2(\R)$.
For example, Kachalov and Kurylev \cite{KK91}
consider inverse scattering  problem, when the potential $V$
satisfies
\[
\lb{4V}
\int_\R(1+|x|)^4|V(x)|dx<\iy.
\]
They prove the recovering problem  using the Gelfand-Levitan equation:
for given S-matrix $S(\l)$ to determine a potential $V$, which satisfies \er{4V} (Theorem 1 and 2 in \cite{KK91}).

Let $V$ satisfy Condition V. Then due to Theorem \er{S1} and \er{HFD}
the S-matrix $S(\l)$ is given by
\[
\label{S01}
S(\l)=S(0)S_1(\l), \qqq S_1(\l)=e^{-2i\l \f_{sc}'(0) }\lim_{r\to +\infty}\prod_{|\l_n|\leq
r}\left({1-{\l\/\ol\l_n}\/1-{\l\/\l_n}}\right) e^{{\l\/\ol\l_n}-{\l\/\l_n}}, \quad \l\in  \C,
\]
uniformly on any compact subset of ${\C}$, where $S(0)={\ol D_+(0)\/D_+(0)}=e^{-2i\f_{sc}(0) }$.
Asymptotics \er{A2} implies $S(\l)\to 1$ as $\l\to \pm\iy$. Then we deduce that there exists a following limit:
\[
\label{S0}
S(0)=S(\l)\ol S_1(\l)=\lim_{\l\to -\iy} \ol S_1(\l).
\]
Thus due to \er{S01} and \er{S0} the S-matrix is expressed in terms of
resonances only.

Now we consider two perturbed Stark operators $H_j=H_0+V_j, j=1,2$  on $L^2(\R)$, where  the potential $V_j$ satisfies Condition V.
We denote the S-matrix for $H_0,H_j$ by $S_j(\l)$.
Assume that $H_1$ and $H_2$ have the same resonances. Then by the above
results, the S-matrices for $H_1$ and $H_2$ coincide, i.e.,
$S_1=S_2$. After this the Kachalov and Kurylev results \cite{KK91}
give $V_1=V_2$.
\BBox

\

\section{Appendix}
\setcounter{equation}{0}

Introduce the Fourier transformation
\[
\lb{}
\wh f(t)={1\/\sqrt{2\pi}}\int_\R f(x)e^{-ixt}dx, \qq t\in \R.
\]

\begin{lemma} \lb{TAY}
Let $V$ satisfy conditions in Theorem \ref{T1} and $
\l\in\ol\C_+$.
Then the following identity and asymptotics
\[
\label{e.3} \Tr Y_0(\l)={e^{i{\pi\/4}}\/ (4\pi
)^{1/2}}\int_{\R}V(x)dx\int_0^\iy e^{-it^3/12}
e^{it(\l-x)}{dt\/t^{1/2}}= {e^{i{\pi\/4}}\/ \sqrt 2}\int_0^\iy
e^{it\l-{i\/12}t^3}\wh V(t) {dt\/t^{1/2}},
\]
\[
\lb{TrY} \Tr Y_0(\l)=i{V_0\/2\sqrt{\l}}+{O(1)\/\l}\qqq as \ |\l|\to
\iy,
\]
hold true, uniformly with respect to $\arg \l\in [0,\pi]$.
\end{lemma}

 {\bf Proof.}  We show \er{e.3}. In view of \er{GSY} we have
$Y_0(\l)\in \cB_1$ for all $\l\in \C_+$. The identity \er{E3} for
$\l\in \C_+$ with $c:={e^{i{\pi\/4}}\/ \sqrt 2}$ gives
$$
\Tr Y_0(\l)=\int_\R V(x)R_0(x,x,\l)dx={c\/
\sqrt{2\pi}}\int_{\R}V(x)dx\int_0^\iy
e^{it(\l-x)-i{t^3\/12}}{dt\/t^{1\/2}} =c\int_0^\iy
e^{it\l-{i\/12}t^3}\wh V(t) {dt\/t^{1\/2}}.
$$

We characterize the properties of  $V$ in terms  of its Fourier
transform $\wh V$:
\[
\lb{FV}
\begin{aligned}
& \wh V, (\wh V)', (\wh V)''  \in L^\iy(\R)\cap L^2(\R),\\
& \wh V(t)={V(0)-e^{-it\g}V(\g)\/it\sqrt{2\pi}}+{1\/it\sqrt{2\pi}}\int_0^\g e^{-itx}V'(x)dx={O(1)\/t}
\qq \as \ t\to \pm \iy.
\end{aligned}
\]

Let us show \er{TrY}. Define cut-off functions $w_0, w_1\in
C^\iy(\R_+)$ by
$$
w_0(t)=\ca 1 & t\in [0,1]\\ 0& t>2\ac,\qqq w_1=1-w_0.
$$
Thus \er{e.3}
 yields the decomposition for $\l\in \ol\C_+$:
\[
\begin{aligned}
& \Tr Y_0(\l)=c\int_0^\iy e^{it\l-{i\/12}t^3}\wh V(t) {dt\/t^{1/2}}=cI_0(\l)+cI_1(\l),\qqq c={e^{i{\pi/4}}\/\sqrt 2},
\\
 & {\rm where} \ I_j=\int_0^\iy w_j(t) e^{i\f(t)}\wh V(t)
{dt\/t^{1/2}},\qq j=0,1,\qq \f(t)=t\l-{t^3\/12}.
\end{aligned}
\]
Let us first consider $I_0$.
Recall that $\wh V(t)$ and has the decomposition:
\[
\begin{aligned}
e^{-{i\/12}t^3}\wh V(t)=\wh V(0)+ t\b(t)\qq \as\ t\to 0,
\end{aligned}
\]
where $\wh V(0)={V_0\/\sqrt{2\pi}}$ and the function $\b$ is entire.
Substituting this asymptotics
into $I_0$ we obtain
\[
\lb{Ix1}
\begin{aligned}
I_0=\int_0^2 w_0(t)e^{it\l-{i\/12}t^3}\wh V(t)
{dt\/t^{1/2}}=I_{01}+I_{02},\qq I_{01}=\wh
V(0)\int_0^2 w_0(t)e^{it\l}{dt\/t^{1/2}}.
\end{aligned}
\]
The stationary phase method gives
\[
\lb{Ix2} I_{01}={V_0\/\sqrt{2\pi}}\int_0^2 w_0(t)
e^{it\l}{dt\/t^{1/2}}={e^{i{\pi/4}}V_0\/\sqrt {2\pi}}\sqrt{\pi\/\l}+{O(1)\/\l^{3\/2}}={cV_0\/\sqrt\l}+{O(1)\/\l^{3\/2}},
\]
and  an integration by parts further  implies
$$
I_{02}=\int_0^2 w_0 e^{it\l}\rt(e^{-{i\/12}t^3}\wh V(t)-\wh V(0)\rt)
{dt\/t^{1/2}}={O(1)\/\l}.
$$
  Thus, if we assume that  $I_1=O(1/\l)$, then   \er{Ix1} and \er{Ix2} yield \er{TrY}.
Hence, it remains to show that  $I_1=O(1/\l)$. To this end we
distinguish  three cases.

 $\bu$ First, let $\l=\m+i\n, \n\ge K|\m|$ for some $K>0$ and let $\n\to \iy$. We have
\[
\lb{I13} |I_1|\le \int_1^\iy e^{-t\n}|\wh V(t)| {dt}|\le \|\wh
V\|_\iy\int_1^\iy e^{-t\n} dt={\|\wh V\|_\iy\/\n} e^{-\n}.
\]

$\bu$  Second, let $\Re \l\to -\iy$. We put $f=w_1(1)\wh V(t)/t^{1/2} $ and here $\f'(t)=\l-{t^2\/4}$. Then due to \er{FV} we have
\[
\lb{Ixx1}
\begin{aligned}
& I_1=\int_1^\iy e^{i\f(t)}f(t)dt= \int_1^\iy
{f(t)\/i\f'(t)} de^{i\f(t)}=i\int_1^\iy
e^{i\f(t)}\rt( {f(t)\/\f'(t)}\rt)'dt={O(1)\/\l}
\end{aligned}
\]
since \er{FV} yields $f', f=O(t^{-3/2}), \ t\to \iy$ and we have the
following estimate
$$
 \rt| \rt( {f(t)\/i(\l-{t^2\/4})}\rt)'  \rt|\le {C\/|\l|\ t^{3\/2}}
\qqq \forall \ t\ge 1,
$$
for some constant $C$

$\bu$  Third, let $\l={k^2+i\n\/4}\in \ol\C_+, \ k\to +\iy$ and $K k^2\ge \n$. Define the cut-off functions $g_0, g_1\in C_0^\iy(\R)$ by
$$
g(s)=\ca 1 & s\in (-1/2,1/2)\\ 0 & s\in \R_+\sm (-1,1) \ac,\qqq
g_1=1-g(t-k).
$$
Then we have
$$
\begin{aligned}
I_1=\int_1^\iy e^{it\l-{i\/12}t^3}f(t)dt=I_{11}+I_{12},\qq I_{11}=\int_1^\iy
e^{it\l-{i\/12}t^3}g(t-k)f(t)dt,\qqq
f=\wh
V(t){w_1(t)\/t^{1\/2}}.
\end{aligned}
$$

Let $F(t)=g_1(t)f(t)$. Similar to \er{Ixx1} we obtain
\[
\lb{J1xw}
\begin{aligned}
I_{12}=\int_1^\iy e^{i\f(t)}F(t)dt= \int_1^\iy
{F(t)\/i\f'(t)} de^{i\f(t)}= -\int_1^\iy
e^{i\f(t)}\rt( {F(t)\/i\f'(t)}\rt)'dt\\
=-\int_1^\iy e^{i\f(t)}\rt({1\/\f'(t)} \rt({F(t)\/\f'(t)}\rt)'\rt)'dt={O(1)\/|\l|}
\end{aligned}
\]
since due to \er{FV} we have
$$
\rt| \rt({1\/\f'(t)} \rt({F(t)\/\f'(t)}\rt)'\rt)'  \rt|\le
{C\/|\l|t^{3\/2}}\qqq \forall \ t\ge 1.
$$

Let us now consider the main term $I_{11}$.
We have for $\l={k^2+i\n\/4}$ and $t=k+s$:
$$
\begin{aligned}
t\l-{i\/12}t^3={1\/12}\rt(3(k+s)(k^2+ i\n)-(k+s)^3\rt)=
{1\/12}\rt(3(k+s)k^2+i3t\n-(k+s)^3\rt)\\
={1\/12}\rt(2k^3+i3t\n-3s^2k-s^3\rt)={k^3\/6}+i\n{k+s\/4}+
-{k\/4}s^2-{s^3\/12}.
\end{aligned}
$$
Define $\p(s,\l)=e^{-i{s^3\/12}}e^{-{\n\/4}(k+s)}g(s)f(k+s)$ and note that $f(k+s)={O(k^{-{3\/2}})}$ as $k\to \iy$, uniformly in $ s\in [-1,1]$.
Due to \er{FV} the stationary phase method gives
\[
\lb{J1xwz}
\begin{aligned}
I_{11}=\int_{k-1}^{k+1} e^{it\l-{i\/12}t^3}g(t-k)f(t)dt
=e^{{i\/6}k^3} \int_{-1}^{1}e^{-{i\/4}ks^2}\p(s,\l)ds
={O(1)\/k^2}={O(1)\/|\l|}.
\end{aligned}
\]
 Combining \er{I13}-\er{J1xwz} we obtain $I_1=O(1/\l)$ and  \er{TrY}.
 \BBox

\bigskip

\footnotesize

\no {\bf Acknowledgments.} \footnotesize Various parts of this paper were
written at the Mathematical Institute of Potsdam  University. E.K. is grateful to the institute for the kind  hospitality. He
is also grateful to Alexei Alexandrov (St. Petersburg), Vladimir
Peller (Michigan) and  Michail Sodin (Tel-Aviv) for stimulating
discussions and useful comments about entire functions. Moreover, he
is also grateful to Oliver Matte (Aarhus) for reading the manuscript
and useful remarks.
The author would like to thank the referee for
useful comments.
His study was supported by the RSF grant  No.
15-11-30007.

\end{document}